\documentclass[12pt]{amsart}
\usepackage{amscd,times,fullpage}

\usepackage{amsthm,amsmath}
\usepackage{amssymb}
\usepackage{graphicx}
\usepackage{amsthm}
\usepackage{enumerate}
\usepackage{tikz-cd}
\usepackage[hidelinks]{hyperref}
\usetikzlibrary{cd}
\usepackage{cleveref}
\usetikzlibrary{positioning}
\usetikzlibrary{matrix}

\newcommand{\R}{\mathbb{R}} 
                   
\newcommand{\Z}{\mathbb{Z}}

\newcommand{\PP}{\mathbb{P}}
\newcommand{\CP}{\mathbb{C}\mathrm{P}}
\newcommand{\N}{\mathbb{N}}

\newcommand{\aut}{\mathrm{Aut}}

\newcommand{\Scal}{\mathrm{Scal}}

\newcommand{\C}{\mathbb{C}}            
\newcommand{\de}{\partial}

\newcommand{\K}{K\"{a}hler }

\newcommand{\OO}{\mathcal{O}}
\newcommand{\BB}{\mathcal{B}_{k+i}}

\newcommand{\ov}[1]{\overline{#1}}

\newcommand{\deb}{\ov\partial}

\newcommand{\di}{{\operatorname{d}}}

\newcommand{\vol}{{\operatorname{vol}}}

\newcommand{\Id}{\operatorname{Id}}
\newcommand{\U}{\operatorname{U}}

\newcommand{\lra}{\longrightarrow}

\newcommand{\D}{\mathcal{D}}

\newcommand{\TM}{\text{T} M}

\newtheorem{theor}{Theorem}
\newtheorem*{theorem}{Theorem}
\newtheorem{prop}{Proposition}
\newtheorem{defin}[prop]{Definition}
\newtheorem{lem}[prop]{Lemma}

\newtheorem{ex}[prop]{Example}
\newtheorem{remark}[prop]{Remark}

\makeatletter
\newtheorem*{rep@theorem}{\rep@title}
\newcommand{\newreptheor}[2]{%
	\newenvironment{rep#1}[1]{%
		\def\rep@title{#2 \ref{##1}}%
		\begin{rep@theorem}}%
		{\end{rep@theorem}}}
\makeatother

\newreptheor{theor}{Theorem}

\usepackage{nicefrac}

\begin{document}

	\title[Any Sasakian structure is approximated by embeddings into spheres]{Any Sasakian structure is approximated by embeddings into spheres}
	
	\author{A.~Loi}
	\address{Dipartimento di Matematica e Informatica, Università degli Studi di Cagliari, Via Ospedale 72, 09124 Cagliari, Italy}
	\email{loi@unica.it}
	\author{G.~Placini}
	\address{Dipartimento di Matematica e Informatica, Università degli Studi di Cagliari, Via Ospedale 72, 09124 Cagliari, Italy}
	\email{giovanni.placini@unica.it}
	
	\subjclass[2010]{53C25; 53C42 ;53C20; 57R18}  
	\keywords{Sasakian geometry; Metric approximation, K\"ahler orbifold embeddings, Sasakian weighted spheres}

	\begin{abstract}
		We show that, for any given $q\geq 0$, any Sasakian structure on a closed manifold $M$ is approximated in the $C^{q}$-norm by structures induced by CR embeddings into weighted Sasakian spheres. In order to obtain this result, we also strengthen the approximation of an orbifold \K form by projectively induced ones given in \cite{rossthomas11b} in the $C^0$-norm to a $C^{q}$-approximation. 
	\end{abstract}
	
	\maketitle
	
	\section{Introduction and statements of the main results}\label{sectionint}
	
	Sasakian geometry is an established field of research in geometry and physics which has received ever growing interest in the past decades.
	The focus on Sasakian geometry is partly justified by its connection to the AdS/CFT theory from the physical point of view.
	From a geometric standpoint the interests is also justified by the abundance of structures underlying a Sasakian manifold and the connection to \K geometry.
	In particular, as \K geometry lies in the intersection of complex, symplectic and Riemannian geometry, a Sasakian structure requires the compatibility of CR, contact and Riemannian structures.
	This allows several different approaches to the study of Sasakian geometry.
	In fact, a Sasakian manifold comes equipped with a characteristic vector field (the Reeb vector field) and a transverse \K geometry.
	When the leaves of the Reeb vector field are compact, its orbit space is a \K orbifold. 
	A fruitful approach is then to study the \K base to obtain informations on the Sasakian manifold itself, see for instance \cite{bandecappellettimontanoloi20,cappellettimontanoloi19,kotschickplacini22,placini21}.
	Other techniques involve the study of the \K cone over a Sasakian manifold. These methods do not require the compactness of the orbits of the Reeb vector field and have proved to be equally effective, see e.g. \cite{collins18,futakionowang09,martellisparksyau08}.
	In this paper we take advantage of both approaches and their interplay.
	
	Given their close relation, many interesting problems in Sasakian geometry originate by analogy from its older sister, \K geometry.
	This is the case in our paper. 
	Namely, we are motivated by a classical result of Tian \cite{tian90} on $C^2$-approximations of polarized \K metrics, later improved by Ruan and Zelditch \cite{ruan98,zelditch98} to a smooth approximation, and its extension due to Ross and Thomas \cite{rossthomas11b} to the orbifold case. In particular, in \cite{rossthomas11b} it is proven that any \K form on a polarized \K orbifold with cyclic quotient singularities can be continuously approximated by embeddings into weighted projective spaces. 
	It is known that not every polarized \K manifold admits an isometric embedding into a complex projective space $\CP^N$.

	In this paper we are interested in the analogous problem in Sasakian geometry. 
	Namely, is there a model space such that all Sasakian structures are approximations of structures induced by embeddings in such a model space? 
	The natural candidate for such a space is a (weighted) Sasakian sphere as it is the Sasakian manifold whose \K base is a (weighted) complex projective space. Such Sasakian structures on spheres are obtained as simple deformations of the standard Sasakian structure, cf. \Cref{ExWeightedSphere}.
	Our main theorem is the following $C^{q}$-approximation theorem.
	
	\begin{theor}\label{TheoMain1}
		Let $(M, \eta, g)$ be a compact Sasakian manifold and fix $q\in\N$. Then there exist a sequence of Sasakian structures induced by CR embeddings $\varphi_k:M\lra S^{2N_k+1}$ into weighted Sasakian spheres which converge to $(\eta,g)$ in the $C^{q}$-norm.
	\end{theor}
	
	Let us briefly comment this result.
	In analogy with the \K setting, the dimension $2N_k+1$ of the model space tends to infinity as $k\lra\infty$.
	The embeddings of \Cref{TheoMain1} do not preserve the Reeb vector field. In fact, the Reeb vector field is the conjugate to the radial direction on the \K cone over $M$. This is determined by the hermitian metric on the cone (see \Cref{EqReebField}) so it is not preserved unless the embedding is in fact a Sasakian embedding.
	Notice that there is no assumption on the regularity of the Sasakian structure in \Cref{TheoMain1}. In fact, there is no assumption other than compactness of the Sasakian manifold.
	
	The main ingredients in the proof of Theorem~\ref{TheoMain1} are:
	\begin{enumerate}[(1)]
		\item $C^\infty$-approximation of irregular structures by quasi-regular ones due to Rukimbira \cite{rukimbira95}.\label{Ingredient1}
		\item $C^{q}$-approximation of a quasi-regular Sasakian structure by Sasakian structures induced by embeddings into weighted Sasakian spheres.\label{Ingredient2}
	\end{enumerate}

	What allows us to get the convergence is the fact that the \K cone does not vary with the structures in point \eqref{Ingredient1}. In fact, Theorem~\ref{TheoMain1} can be rephrased as a statement on embeddings of the \K cone of $M$ into $\C^{N_k+1}\setminus\{0\}$.
	Moreover, the structures induced via embedding are simple deformations (see Definition~\ref{DefinSimpleTransform} below) of the original one.
	All these deformations do not change the CR structure so we get a sequence of CR embeddings into weighted Sasakian spheres that approximate a given structure.
	This is in analogy with the \K setting where the embeddings are holomorphic, since the CR structure can be regarded as the transverse holomorphic structure.

When the structure $(M, \eta, g)$ is quasi-regular, the space of leaves is a polarized \K orbifold $X$. This allows us to work on the \K geometry of the base $X$ and approximate it by embeddings into projective spaces. 
Such approximations given by the asymptotic of the Bergman kernel for the Kodaira-Baily embedding were studied by Dia, Liu and Ma in \cite{DaiLiuMa06}. 
Unfortunately, this does not suit our needs as the Sasakian structures so obtained would be regular. 
Thus we require an \textit{orbifold} embedding (cf. \cite[Section~2.1]{rossthomas11b}) in a \textit{weighted} projective space.
This was done by Ross and Thomas \cite{rossthomas11b} in the $C^2$-norm. The main technical difficulty of part \eqref{Ingredient2} is to extend their result to a $C^{q}$-convergence. This is done in Theorem~\ref{TheoMain2}.

When the Sasakian structure on $M$ is regular, that is, the Reeb action is free, it is natural to ask whether one can get a similar result under the requirement that the model space is a standard Sasakian sphere.
We address this problem, both in the compact and noncompact setting, in a related work \cite{placini22}.

\subsection*{Related Works}	
Recently, a closely related problem was considered by Herrmann, Hsiao, Marinescu, and Shen in \cite{HerrmannHsiaoMarinescuShen22}. There, using different methods, the authors provide several results on CR embeddings into spheres with the weighted CR structure. Remarkably, they obtain smooth approximation of CR structures. In particular, one should compare our \Cref{TheoMain1} with \cite[Theorem~6.5]{HerrmannHsiaoMarinescuShen22}
One should compare \Cref{TheoMain1} also with the work of Ornea and Verbitsky \cite{orneaverbitsky05,orneaverbitsky07a,orneaverbitsky07b,orneaverbitsky22} dealing with CR embeddings of Sasakian manifolds into spheres. In particular, our result can be regarded as a strengthening of \cite[Theorem~11.21]{orneaverbitsky22}.
	More generally, CR embeddings $\Phi_m:M\lra \C^{N_m}$ of any CR manifold $M$ were produced in \cite{HerrmannHsiaoLi17,HerrmannHsiaoLi22} for all  $m$ sufficiently large. Moreover, the authors studied the asymptotic of the Szeg\H{o} kernel. However, this cannot be applied to the Sasakian case in a straightforward manner as the Szeg\H{o} kernel considered in \cite{HerrmannHsiaoLi17,HerrmannHsiaoLi22} relates to the pullback of the standard metric of $\C^{N_m}$ rather than the cone metric inducing the Sasakian structure on the weighted sphere.

\subsection*{Organization of the paper}	The remainder of the paper is organized as follows.
	In Section~\ref{sectionbackground} we review Sasakian geometry and define the basic objects needed in the remainder of the article.
	In Section~\ref{SectionAsymptotics} we recall some notions on weighted Bergman kernels and their asymptotic expansions from \cite{rossthomas11b} and use them to prove Theorem~\ref{TheoMain2} in Section~\ref{sectionprooforbifold}.
	Finally Section~\ref{sectionproof1} is dedicated to the proof of Theorem~\ref{TheoMain1}.

	\vspace{3mm}
	\noindent We would to thank Julius Ross for his interest in our work and his valuable comments.

	\section{Sasaki manifolds and \K cones}\label{sectionbackground}

	Sasakian geometry can be presented in terms of contact metric geometry as in the monograph of Boyer and Galicki \cite{boyergalicki08} or via the associated \K cone, see e.g. \cite{hesun16,martellisparksyau08,sparks11} for easily readable references. We will present both formulations because a combination of the two is most suitable to our purpose.
	In the following all manifolds and orbifolds are assumed to be connected, closed and oriented.

	A \textit{K-contact structure} $(\eta,\Phi,R,g)$ on a manifold $M$ consists of a contact form $\eta$ and an endomorphism 
	$\Phi$ of the tangent bundle $\TM$ satisfying the following properties:
	\begin{enumerate}
		\item[$\bullet$] $\Phi^2=-\Id+R\otimes\eta$ where $R$ is the Reeb vector field of $\eta$,
		\item[$\bullet$] $\Phi_{\vert\D}$ is an almost complex structure compatible with the symplectic form $\di\eta$ on $\D=\ker\eta$,
		\item[$\bullet$] the Reeb vector field $R$ is Killing with respect to the metric $g=\frac{1}{2}\di\eta\circ\Id\otimes\Phi+\eta\otimes\eta$.
	\end{enumerate}   
	Given such a structure one can consider the almost complex structure $J$ on the Riemannian cone $\big( M\times\R^+,t^2g+\di t^2\big)$ given by
	\begin{enumerate}
		\item[$\bullet$] $J=\Phi$ on $\D=\ker\eta$, and
		\item[$\bullet$] $R=J(t\partial_t)_{\vert_{\{t=1\}}}$.
	\end{enumerate} 
	A \textit{Sasakian structure} is a K-contact structure $(\eta,\Phi,R,g)$ such that the associated almost complex structure $J$ is integrable, and therefore $\left( M\times\R^+,t^2g+\di t^2,J\right)$ is K\"ahler. 
	A \textit{Sasakian manifold} is a manifold $M$  equipped with a Sasakian structure $(\eta,\Phi,R,g)$.
	
	Equivalently, one can define Sasakian manifolds in terms of \K cones. 
	Namely, a Sasakian structure on a compact smooth manifold $M$ is defined to be a \K cone structure on $Y:=M\times\R^+$.
	That is, A \K structure $(g_Y,J)$ on $Y$ of the form 
	$g_Y=t^2g+\di t^2$ where $t$ is the coordinate on $\R^+$ and $g$ a metric on $M$.
	Then $(Y,g_Y,J)$ is called \textit{the \K cone of $M$} which, in turn, is identified with the submanifold $\{t=1\}$. The \K form on $Y$ is then given by
	$$\Omega_Y=\dfrac{i}{2}\de\deb t^2.$$
	The Reeb vector field on $Y$ is defined as 
	$$R=J(t\de_t).$$
	This defines a holomorphic Killing vector field on $Y$ with metric dual $1$-form
	$$\eta=\dfrac{g_Y(R,\cdot)}{t^2}=\di^c\log t=i(\deb-\de)\log t$$
	where $d^c=i(\deb-\de)$.
	With a slight abuse of notation we write $R$ and $\eta$ to indicate both the objects on $Y$ and their restrictions to $M$.
	Notice that $J$ induces an endomorphism $\Phi$ of $\TM$ by setting 
	\begin{enumerate}
		\item[$\bullet$] $\Phi=J$ on $\D=\ker\eta_{\vert_{\TM}}$, and
		\item[$\bullet$] $\Phi(R)=0$.
	\end{enumerate} 
	Equivalently, the endomorphism $\Phi$ is determined by $g$ and $\eta$ by setting
	$$g(X,Z)=\dfrac{1}{2}\di\eta(X,\Phi Z)\ \ \mbox{ for }\ X,Z\in\D.$$
	It is easy to see that, when restricted to $M=\{t=1\}$, $(\eta,\Phi,R,g)$ is a Sasakian structure in the contact metric sense whose \K cone is $(Y,g_Y,J)$ itself. 
	
	The Reeb vector field $R$ defines a foliation on $M$, called the Reeb foliation.
	A very important dichotomy of Sasakian manifolds $M$ is given by compacteness of the leaves of the Reeb foliation on $M$.
	Namely, when the leaves of the Reeb vector field are noncompact, the Sasakian structure is called \textit{irregular}. If the orbits of $R$ are compact, there is a further distinction. That is, the Sasakian structure is called \textit{regular} if all the orbits have the same period and \textit{quasi-regular} otherwise.
	
	Since $g$ and $\eta$ are invariant for $R$ the Reeb foliation is transversally \K in the sense that the distribution $\D$ admits a \K structure $(g^T,\omega^T,J^T)$ which is invariant under $R$.
	Explicitly, we have 
	$$\omega^T=\dfrac{1}{2}\di\eta,\ \ J^T=\Phi_{\vert_{\D}}\ \mbox{ and } \ g^T(X,Z)=\dfrac{1}{2}\di\eta(X,J^T Z)=g_{\vert_{\D}}.$$
	In particular, we have
	\begin{equation}\label{EqTransvFormIsCurv}
		\omega^T=\dfrac{1}{2}\di\eta=\dfrac{i}{2}\di(\deb-\de)\log t=\dfrac{i}{2}\de\deb\log t^2.
	\end{equation}
	
	Regular and quasi-regular Sasakian manifold are fairly well understood due to the following result 
	
	\begin{theorem}[\cite{boyergalicki08}]\label{TheorStructure}
		Let $(M,\eta,\Phi,R,g)$ be a quasi-regular  compact Sasakian manifold. Then the space of leaves of the Reeb foliation $(X,\omega)$ is a compact \K cyclic orbifold with integral \K form $\frac{1}{\pi}\omega$ 
		so that the projection $\pi:M\lra X$ is a Riemannian submersion. Moreover, $X$ is a smooth manifold if and only if the Sasakian structure on $M$ is regular. 
		
		Viceversa, any principal $S^1$-orbibundle $M$ with Euler class $-\frac{1}{\pi}[\omega]\in H^2_{orb}(X,\Z)$ over a compact \K cyclic orbifold $(X,\omega)$ admits a Sasakian structure.
	\end{theorem}
	We refer to \cite{boyergalicki08} for the basics of \K orbifolds too.
	This result allows us to reformulate the geometry of a quasi-regular Sasakian manifold $M$ in terms of the  \K geometry of the space of leaves $X$. 
	We will illustrate in detail this correspondence for its importance in the remainder of the paper.
	In order to do so, let us first introduce the concept of $\D$-homothetic transformation of a Sasakian structure.
	
	\begin{defin}[$\D$-homothety or transverse homothety {\cite[Section~7.3]{boyergalicki08}}]
		Let $(M,\eta,\Phi,R,g)$ be a Sasakian manifold and $a\in\R$ a positive number. One can define the Sasakian structure $(\eta_a,\Phi_a, R_a,g_a)$ from $(\eta,\Phi,R,g)$  as
		$$ \eta_a=a\eta,\ \ \ \Phi_a=\Phi,\ \ \ R_a=\dfrac{R}{a},\ \ \  g_a=ag+(a^2-a)\eta\otimes\eta=ag^T+\eta_a\otimes\eta_a.$$
	\end{defin}
	Equivalently, the Sasakian structure $(\eta_a,\Phi_a, R_a,g_a)$ on $M$ can be obtained from the \K cone by setting the new coordinate $\widetilde{t}=t^a$.

	Now let the quasi-regular Sasakian manifold $(M,\eta,\Phi,R,g)$ be given and consider the projection $\pi:(M,g)\lra(X,\omega)$ given above. 
	Notice that $\pi$ locally identifies the contact distribution $\D$ with the tangent space of $X$. Therefore,  up to $\D$-homothety, we have that $\pi^*(\omega)=\frac{1}{2}\di\eta$. 
	Moreover, the endomorphism $\Phi$ determines the complex structure on $X$ and $g$ induces the \K metric $g_\omega$ compatible with $\omega$, i.e. $g^T=\pi^*g_\omega$.
	
	In this case the class $\frac{1}{\pi}[\omega]\in H^2_{orb}(X,\Z)$ defines a orbiample line bundle over $X$ in the sense of \cite[Definition~2.7]{rossthomas11b}.
	Moreover, the cone $Y=M\times\R^+$ is identified with $L^*$ without the zero section in the following way. Let $h$ be a hermitian metric on $L$ such that
	\begin{equation}\label{EqCurvForm}
		\omega=-\dfrac{i}{2}\de\deb\log h.
	\end{equation}
	Then its dual $h^*$ on $L^*$ defines the second coordinate $(p,t)\in M\times\R^+=L^*\setminus\{0\}$ by
	\begin{align}\label{EqReebField}
		t:L^*\setminus\{0\} &\lra \R^+\\ 
	\nonumber	(p,v)&\mapsto\vert v\vert_{h^*_p}
	\end{align}
	where $v$ is a vector of $L^*$ in the fiber over $p$.
	With this notation the \K form on the \K cone $\left( M\times\R^+,t^2g+\di t^2,J\right)$ is given by
	
	\begin{equation}
		\Omega=\dfrac{i}{2}\de\deb t^2.
	\end{equation}
	The Sasakian structure can be read from this data as
	\begin{equation}\label{EqFromHermitianToSasakian}
		\omega^T=-\dfrac{i}{2}\de\deb\log h,\ \ \ \ \eta=i(\deb-\de)\log t.
	\end{equation}

	Therefore, the choice of a hermitian metric $h$ on a orbiample line bundle $L$ over a compact \K orbifold $X$ completely determines a Sasakian structure on the $U(1)$-orbibundle associated to $L^*$. Further, this is a smooth manifold if the local uniformizing groups of the orbifold inject into $U(1)$.
	The Sasakian manifold so obtained is called a \textit{Boothby-Wang bundle} over $(X,\omega)$ (cf. for instance \cite[Chapter~7]{boyergalicki08}).

	The most basic example is the standard Sasakian structure on $S^{2n+1}$, that is, the Boothby-Wang bundle determined by the Fubini-Study metric $h=h_{FS}$ on $\OO(1)$ over $\CP^n$. We give the details of this construction to further illustrate the formulation above.
	\begin{ex}[Standard Sasakian sphere]\label{ExStandardSphere}
		Let $h=h_{FS}$ be the Fubini-Study hermitian metric on the holomorphic line bundle $\OO(1)$ over $\CP^n$. Recall that its dual metric $h^*$ on $\OO(-1)\setminus\{0\}=\C^{n+1}\setminus\{0\}$ is given by the euclidean norm.
		This defines a coordinate $t$ on the \K cone $\OO(-1)\setminus\{0\}=\C^{n+1}\setminus\{0\}=S^{2n+1}\times \R^+$. Namely, for coordinates $z=(z_0,z_1,\ldots,z_n)$ on $\C^{n+1}$ we have
		\begin{align*}
			t:\C^{n+1} &\lra \R^+\\ 
			z&\mapsto\vert z\vert=\sqrt{\sum_{i=0}^{n}z_i\overline{z}_i}
		\end{align*}
		Now the \K metric on the cone is nothing but the flat metric
		$$	\Omega_{flat}=\dfrac{i}{2}\de\deb t^2=\dfrac{i}{2}\sum\di z_i\wedge\di \overline{z}_i.$$
		The Reeb vector field $R_0$ and the contact form $\eta_0$ read
		$$	R_0=J(t\de_t)=i\sum z_i\de_{z_i}-\overline{z}_i\de_{\overline{z}_i},\ \ \ \eta_0=i(\deb-\de)\log t=\dfrac{i}{2t^2}\sum z_i\di\overline{z}_i-\overline{z}_i\di z_i.$$
		It is clear that, when restricted to $S^{2n+1}$, $\eta_0$ and $R_0$, together with the round metric $g_0$ and the restriction $\Phi_0$ of $J$ to $\ker\eta_0$ give a Sasakian structure on $S^{2n+1}$.
		This corresponds exactly to the Hopf bundle $S^{2n+1}\lra\CP^n$. Moreover, we have
		$$\pi^*\omega_{FS}=\omega^T=\frac{1}{2}\di\eta_0=\dfrac{i}{2\vert z\vert^4}\sum_i\vert z_i\vert^2\di z_i\wedge\di \overline{z}_i-\sum_{i,j}\overline{z}_i z_j\di z_i\wedge\di \overline{z}_j$$
		where $\pi:\C^{n+1}\setminus\{0\}\lra\CP^n$ is the standard projection and the Fubini-Study form is normalized to give $\vol(\CP^n)=\pi^n$.
	\end{ex}
	\subsection{Some Sasakian deformations and weighted Sasakian spheres}
	We begin this section by recalling some well known classes of deformations of Sasakian structures.
	Let us start with transverse \K transformations.
	Namely, given a \K cone $Y=M\times\R^+$ we consider all \K metrics on $(Y,J)$ that are compatible with the Reeb field $R$. In other terms, these are potentials $\widetilde{t}^2$ such that $t\de_t=\widetilde{t}\de_{\widetilde{t}}$.
	This means that the corresponding \K and contact forms satisfy 
	$$\widetilde{\Omega}=\Omega+i\de\deb e^{2f},\ \ \ \widetilde{\eta}=\eta+d^cf$$
	for a function $f$ invariant under $\de_t$ and $R$. Such functions are called \textit{basic functions}.
	We still need to identify the manifolds $\{\widetilde{t}=1\}$ and $\{t=1\}$.
	This is done by means of the diffeomorphism
	\begin{align*}
		F:Y&\lra Y\\
		(p,t)&\mapsto\left(p,te^{-f(p)}\right)
	\end{align*}
	which maps $\{t=1\}$ to $\{t=e^{-f(p)}\}=\{\widetilde{t}=1\}$. It is elementary to check that $\eta,\ R$ and $d^cf$ are invariant under $F$ so that $\widetilde{\eta}=\eta+d^cf$ holds on $M$.
	Furthermore, the transverse \K forms are related by $\widetilde{\omega}^T=\omega^T+i\de\deb f$.
	Notice that when the Sasaki structure is quasi-regular basic functions correspond to function on the base orbifold $X$. Thus, if $t$ comes from a hermitian metric $h^*$ on $L^*$, such a transformation is given by replacing $h^*$ with $e^fh^*$ for a function $f:X\lra\C$ such that $\omega+i\de\deb f>0$. This is equivalent to picking a different \K form $\widetilde{\omega}$ in the same class as $\omega$.
	We summarize the above discussion in the following definition (see e.g. \cite[Proposition~1.4]{sparks11}).
	\begin{defin}[Transverse \K deformations]\label{DefinTransverseDeformation}
		Let  $(M,\eta,R,g,\Phi)$ be a Sasakian manifold with \K cone $(Y,J)$ and \K potential $t^2$. A transverse \K transformation is given by replacing $t$ with $\widetilde{t}=e^{f}t$ for a basic function $f$ and leaving $(Y,J,R)$ unchanged. When the Sasaki structure is quasi-regular and given as in \eqref{EqFromHermitianToSasakian}, a transverse \K transformation is given by replacing $h^*$ with $e^fh^*$.
	\end{defin}
	We now focus on deformations induced by a modification of the Reeb vector field rather than the the \K potential. Such deformations are called Type-I deformations.
	Notice that a transformation of the Sasakian structure $(M,\eta,R,g,\Phi)$ naturally induces an automorphism of the \K cone $(Y,J)$. 
	Therefore, the group $\aut(M,\eta,R,g,\Phi)$ of diffeomorphisms of $M$ preserving the Sasaki structure is a subgroup of $\aut(Y,J)$.
	Fix a maximal torus $\mathbb{T}$ in $\aut(M,\eta,R,g,\Phi)$ and pick an element $R'$ of its Lee algebra $\mathfrak{t}$ such that $\eta(R')>0$. 
	Then there exists a $\mathbb{T}$-invariant coordinate on the \K cone $Y$ inducing a $\mathbb{T}$-invariant \K cone metric on $(Y,J)$ with Reeb vector field $R'$, see e.g. \cite[Lemma~2.2]{hesun16}.
	We give the deformed Sasakian structure in the following
	\begin{defin}[Type-I deformations]\label{DefinTypeI}
		Let  $(M,\eta,R,g,\Phi)$ be a Sasakian manifold. A type-I deformation is given by a choice of a maximal torus $\mathbb{T}$ in $\aut(M,\eta,R,g,\Phi)$ and an element $R'$ of its Lee algebra $\mathfrak{t}$ such that $\eta(R')>0$. The deformed structure is then given by
		$$\eta'=\dfrac{\eta}{\eta(R')},\ \ \Phi'=\Phi-\Phi R'\otimes\eta',\ \ g=\eta'\otimes\eta'+\dfrac{1}{2}\di\eta'(\Id\otimes\Phi').$$
	\end{defin}
	Notice that if a Sasakian manifold is irregular, then $R$ is an irrational element in the Lie algebra $\mathfrak{t}$ of a maximal torus $\mathbb{T}$ in which it is contained.
	This allows us to reformulate a classical result of Rukimbira \cite{rukimbira95} in terms of the \K cone. A proof of the following proposition can be found, for instance, in \cite[Corollary~2.8]{collins18}.
	\begin{prop}\label{PropRukimbira}
		Let  $(M,\eta,R,g,\Phi)$ be a Sasakian manifold such that its Reeb vector field is irrational in a maximal torus $\mathbb{T}$ in $\aut(M,\eta,R,g,\Phi)$ (i.e. an irregular Sasakian manifold).
		Then there exists a sequence of Reeb fields $R_j\in\mathfrak{t}$ such that the corresponding type-I deformations $(M,\eta_j,R_j,g_j,\Phi_j)$ are quasi-regular and converge smoothly to $(M,\eta,R,g,\Phi)$.
	\end{prop}

	We now conclude the section with the definition of weighted Sasaki sphere.
	
	\begin{defin}[Simple deformation]\label{DefinSimpleTransform}
		Let  $(M,\eta,R,g,\Phi)$ be a Sasakian manifold with \K cone $(Y,J)$ and fix a maximal torus $\mathbb{T}$ in $\aut(M,\eta,R,g,\Phi)$.
		A simple deformation of $(M,\eta,R,g,\Phi)$ is a Sasakian structure on $M$ induced by a \K cone metric on $(Y,J)$ with Reeb vector field $R'\in\mathfrak{t}$ such that $\eta(R')>0$.
	\end{defin}
	
	Notice that a simple deformation amounts to a type-I deformation followed by a transverse \K deformation. Since neither of the two changes the CR structures, simple deformations preserve the underlying CR structure.

	\begin{ex}[Weighted Sasaki sphere]\label{ExWeightedSphere}
		Let $(S^{2n+1},\eta_0,R_0,g_0,\Phi_0)$ be the standard structure coming from the flat \K metric on $\C^{n+1}\setminus\{0\}$.
		Then $\aut(S^{2n+1},\eta_0,R_0,g_0,\Phi_0)=U(n+1)$ and we fix the maximal torus $\mathbb{T}^{n+1}\subset U(n+1)$ acting by diagonal elements.
		The Lie algebra $\mathfrak{t}$ of $\mathbb{T}$ is generated by the elements $R_j=i \left(z_j\de_{z_j}-\overline{z}_j\de_{\overline{z}_j}\right)$.
		Therefore an element $R_w\in\mathfrak{t}$ satisfies the positivity condition $\eta(R_w)>0$ if and only if $R_w=\sum_jw_jR_j$ where $w=(w_0,w_1,\ldots,w_n)$ is a positive vector in $\R^{n+1}$.
		
		A simple Sasakian structure on $S^{2n+1}$ is a Sasakian structure that can be obtained as simple deformation of the standard Sasakian sphere with Reeb field $R_w$ as above. The sphere $S^{2n+1}$ endowed with a simple Sasakian structure is called a \textit{weighted Sasaki sphere}.
		Notice that when the weights $w_j$ are rational the Sasakian structure is quasi-regular and the base orbifold is the quotient of $\C^{n+1}\setminus\{0\}$ by the $\C^*$-action
		$$\lambda\cdot(z_o,z_1,\ldots,z_n)=(\lambda^{w_0}z_0,\lambda^{w_1}z_1,\ldots,.\lambda^{w_n}z_n).$$
		This is a \textit{weighted projective space} $\CP(w)$ and the weighted Sasakian structure corresponds to a hermitian metric $h$ on  $\OO(1)$ such that its curvature pulls back to $\frac{1}{2}\di\eta_w$ on $S^{2n+1}$.
		Notice that up to transverse homothety one can assume that the weights $w_j$ are positive integers.
		Moreover, one gets (a $\D$-homothety of) the standard Sasakian structure if and only if all weights are equal.
	\end{ex}

	This realizes $S^{2n+1}$ as the $U(1)$ bundle associated to the line bundle $\OO(-1)$ over the weighted projective space $\CP(w)$, that is, the quotient of $C^{n+1}\setminus\{0\}$ by the weighted $\C^*$ action above.
	The weighted Sasakian structure is then determined by a metric $h$ on  $\OO(1)$ such that its curvature is a \K form on $\CP(w)$.
	Notice that there is more than one natural choice for $h$, see \cite[Section~3]{rossthomas11b}.
	
	In analogy with the unweighted case, one could wish to restrict to the case where $h$ is the Fubini-Study hermitian metric whose curvature is the Fubini-Study \K metric on $\CP(w)$.
	This is what is done, for instance, in \cite[Example~7.1.12]{boyergalicki08} while our definition allows different choices for $h$.
	Nevertheless, the metric used in our Theorem~\ref{TheoMain1} coincides exactly with the one in \cite{boyergalicki08}.
	In fact, it was noted in \cite[Theorem~1.6]{sparks11} that $\frac{1}{2}t_w^2$ is the Hamiltonian function for the Reeb vector field $R_w$ (where $t_w$ is the coordinate on $\C^{n+1}\setminus\{0\}$ coming from the type-I deformation).
	This implies that \K orbifold $(\CP(w),\omega)$ of the weighted Sasakian sphere is the one obtained as \K reduction of the \K cone $\C^{n+1}\setminus\{0\}$ under the $U(1)$-action of the Reeb field $R_w$ (cf. Definition~\ref{DefFubiniOrbi} below).

	\section{Asymptotic expansion of the weighted Bergman kernel}\label{SectionAsymptotics}
	
	\subsection{Kodaira embedding for \K orbifolds}\label{SubSectionEmbeddings}
	In this section we recall some results from \cite{rossthomas11a,rossthomas11b} for clarity and future reference. In order to improve comprehensibility we follow their notation closely.
	
	Let $(X,L)$ be a polarized cyclic orbifold and let $h$ be a hermitian metric on $L$ with curvature $2\pi\omega$ (cf. Remark~\ref{RemarkConventions} below).
	Consider the weighted vector space
	$$V=\bigoplus_iH^0(L^{k+i})^*$$
	where $i$ ranges from $1$ to the order $m=\mathrm{ord(X)}$ of the orbifold $X$ and the $i$-th summand has weight $k+i$.
	Then, for large enough $k\gg0$ we have an embedding $\phi_k:X\lra \mathbb{P}(V)$ such that $\phi_k^*\OO(1)=L$ (\cite{rossthomas11a}) given by
	$$\phi_k(x)=\left[\oplus_i\mathrm{ev}_x^{k+i}\right]$$
	where $\mathrm{ev}_x^{k+i}$ is the evaluation element at $x$ of $H^0(L^{k+i})^*$ which sends a section $s$ of $L^{k+i}$ to $s(x)\in L^{k+i}$ identified with $\C$ via a local trivialization.
	
	We now explain the relation between the hermitian metric $h$ on $L$ and the Fubini-Study metric on $\PP(V)$. 
	Let $\vert\cdot\vert_V$ be a hermitian metric on $V$ such that the summands $H^0(L^{k+i})^*$ and $H^0(L^{k+j})^*$ are orthogonal for $i\neq j$, that is, a collection of hermitian metrics $\vert\cdot\vert_{k+i}$ on each summand of $V$.
	The group $\U(1)$ acts on each $H^0(L^{k+i})^*$ with weight $k+i$. For an element $v=\oplus_iv_{k+i}\in V$ the moment map of the $\U(1)$ action is given by
	\begin{equation}
		\mu_{\U(1)}(v)=\dfrac{1}{2}\left(\sum_i(k+i)\vert v_{k+i}\vert^2 -c\right)\ \ \mbox{with} \ \ c:=\sum_i(k+i)c_ih_0(L^{k+i})
	\end{equation}
	where $c_i$ are arbitrary positive real constants.
	\begin{defin}\label{DefFubiniOrbi}
		The Fubini-Study metric $\omega_{FS}$ on $\PP(V)$ is defined as $\frac{1}{c}$ times the metric induced by $\vert\cdot\vert_V$ on $\PP(V)=\mu_{\U(1)}^{-1}(0)/\U(1)$ by \K reduction.
	\end{defin}
	
	Moreover, the Fubini-Study metric on $\OO_{\PP(V)}(-1)$ in \cite{rossthomas11b} is defined so that the vectors of $\mu_{\U(1)}^{-1}(0)$ are unitary. Namely, we have the following
	\begin{defin}\label{DefHermitianMetric}
		The Fubini-Study hermitian metric $h_{FS}$ on $\OO_{\PP(V)}(-1)$ is defined as $\vert v\vert_{h_{FS}}:=\lambda(v)^{-1}$ where $\lambda(v)$ is the unique positive real number such that $\lambda(v)\cdot v\in\mu_{\U(1)}^{-1}(0)$.
		In other terms, $\lambda(v)$ is the unique positive real solution to $\sum_i(k+i)\lambda(v)^{2(k+i)}\vert v_{k+i}\vert^2=c$.
	\end{defin}
	By abuse of notation we use $h_{FS}$ for the induced metric on $\OO_{\PP(V)}(i)$ as well.
	Unfortunately the curvature $2\pi\omega_{h_{FS}}:=i\de\deb\log h_{FS}$ of $h_{FS}$ on $\OO_{\PP(V)}(1)$  is not given by $2\pi\omega_{FS}$. In fact, their difference is described in the next
	
	\begin{lem}{{\cite[Lemma~3.6]{rossthomas11b}}}\label{LemFubiniStudy}
		In the notation above we have
		$$\omega_{FS}=\omega_{h_{FS}}+\dfrac{i}{2c}\de\deb f$$
		where $f=\sum_i\sum_\alpha\vert t^i_\alpha\vert^2_{h_{FS}}$
		and  $\{t^i_\alpha\}$ is an orthonormal basis of $(V^*,\vert\cdot\vert_V)$.
	\end{lem}
	\begin{remark}\rm\label{RemDefRealFubini}
		In particular this implies that there exists a hermitian metric $h'_{FS}=e^{\frac{f}{2c}}h_{FS}$ on $\OO_{\PP(V)}(1)$ whose curvature is exactly $2\pi\omega_{FS}$.
	\end{remark}
	
	\begin{remark}\rm\label{RemarkConventions}
		Notice that the convention $2\pi\omega_{h}=i\de\deb\log h$ used here for the curvature form $\omega_h$ of a hermitian metric $h$ on a line bundle $L$ differs from Equation \eqref{EqCurvForm} by a factor of $-\pi$.
		This is of no consequence for the convergence results we are aiming for. Thus we follow the convention of \cite{rossthomas11b} throughout this section and Section~\ref{sectionprooforbifold}.
	\end{remark}
	Notice that there is a correspondence between pairs $(h,\omega)$, consisting of a hermitian metric $h$ on $L$ and an orbifold \K metric $\omega$ on $X$ such that $2\pi\omega$ is the curvature of some hermitian metric on $L$, and metrics $\vert\cdot\vert_V$ on $V=\bigoplus_iH^0(L^{k+i})^*$. 
	Namely, given a pair $(h,\omega)$ as above, one can define a metric on $V$ (which following \cite{rossthomas11b} we denote by $\mathrm{Hilb}(h,\omega)$) by setting
	$$\vert s\vert^2_{\mathrm{Hilb}(h,\omega)}=\dfrac{1}{c_i\vol}\int_{X}\vert s\vert^2_{h}\dfrac{\omega^n}{n!}$$
	where $s$ is a section of $L^{k+1}$ and
	$$\vol:=\int_X\dfrac{c_1(L)^n}{n!}=\int_X\dfrac{\omega^n}{n!}.$$ 
	
	Conversely, a metric $\vert\cdot\vert_V$ on $V$ gives a pair of metrics $(h_{FS},\omega_{FS})$ on $\OO_{\PP(V)}(1)$ and $\PP(V)$ as constructed above. In turn they determine a pair $(h,\omega)=(\phi_k^*h_{FS},\phi_k^*\omega_{FS})$.
	
	If the pair $(h,\omega)$ is fixed and is inducing the metric $\vert\cdot\vert_V=\mathrm{Hilb}(h,\omega)$,  then we denote the pair $(\phi_k^*h_{FS},\phi_k^*\omega_{FS})$ by $(h_{FS,k},\omega_{FS,k})$. From now on we will assume that this is the case.
	
	\subsection{The weighted Bergman kernel}\label{SubSectionBergmanKernel}
	
	Let $(X,L)$ be a polarized orbifold, $h$ a hermitian metric on $L$ and $\omega$ a \K metric on $X$. As in the smooth case one can give the following
	\begin{defin}
		The Bergman kernel of $L^{k+i}$ is defined to be the function $\BB(h,\omega)=\sum_\alpha\vert t_\alpha\vert^2_h$ where $\{t_\alpha\}$ is an orthonormal basis of $H^0(L^{k+i})$ with respect to $\mathrm{Hilb}(h,\omega)$.
	\end{defin}
	
The functions $\BB$ satisfy a crucial property that will allow us the control we need to prove Theorem~\ref{TheoMain2}. This is proven in \cite[Corollary~4.10]{rossthomas11a} and, alternatively, can be deduced from \cite[Theorem~5.4.11]{marinescuma}.
	
	\begin{lem}
		If $D$ be a differential operator of order $p$, then
		\begin{equation}\label{EqDerivBergKer}
			D\mathcal{B}_k=O(k^{n+p})
		\end{equation}
		uniformly on $X$ for large $k$.
	\end{lem}
	
	This is not the only Bergman kernel we are interested in, though. In fact, we consider the weighted Bergman kernel, that is, the version of it related to orbifold embeddings. In order to define the weighted Bergman kernel, let us pick a basis  $\{s^i_\alpha\}$ of $\oplus_iH^0(L^{k+i})$ orthonormal with respect to the $L^2$-metric and a basis  $\{t^i_\alpha\}$ orthonormal for $\mathrm{Hilb}(h,\omega)$. For instance, one could take $t^i_\alpha=\sqrt{c_i\vol}s^i_\alpha$.
	
	\begin{defin}
		The weighted Bergman kernel is defined to be the function 
		
		\begin{equation}\label{EqBergmanL2}
			B_k=B_k(h,\omega)=\vol\sum_ic_i(k+i)\sum_\alpha\vert s^i_\alpha\vert^2_h
		\end{equation}

		or, equivalently,
		\begin{equation}\label{EqBergmanHilb}
			B_k=\sum_i(k+i)\sum_\alpha\vert t^i_\alpha\vert^2_h=\sum_i(k+i)\BB
		\end{equation}
		where $i$ runs from $1$ to $m=\mathrm{ord(X)}$.
	\end{defin}
	
	A crucial element for the proof of Theorem~\ref{TheoMain2} is the asymptotic expansion of the weighted Bergman kernel $B_k=\sum_i(k+i)\BB$ proven in \cite[Theorem~1.7 and Remark~4.13]{rossthomas11a}. Namely, for coefficients $c_i$ given by 
	\begin{equation}\label{EqCi}
		\sum_ic_it^i:=\left(t^{m-1}+t^{m-2}+\cdots+1\right)^p
	\end{equation}
	there is an asymptotic expansion
	\begin{equation}\label{EqExpansionCq}
		B_k=b_0k^{n+1}+b_1k^n+\cdots\ \mbox{for}\ k\rightarrow\infty
	\end{equation}
	for some smooth functions $b_i$. 
	Moreover, for a given $q\in\N$ we can pick the coefficients $c_i$ given as in \Cref{EqCi} with $p>r+q$. This will ensure that the expansions holds in $C^q$ up to order $O(k^{n+1-r})$. Assume from now on that $q\in\N$ is fixed. We will always assume that $p$ is chosen large enough so that the expansions holds in $C^q$.
	
	In our case $2\pi\omega$ is the curvature of $h$ so that \cite[(1.11) and Corollary~4.12]{rossthomas11a} give
	
	\begin{equation}
		b_0=\vol\sum_ic_i,\ \ \ \ \ b_1=\vol\sum_ic_i\left((n+1)i+\dfrac{1}{2}\Scal(\omega)\right)
	\end{equation}
	where $\Scal(\omega)$ is the scalar curvature of $\omega$.
	In particular, $b_0$ is constant over $X$. 
	Moreover, integrating over $X$ one sees that $c=\sum_ic_i(k+i)h^0(L^{k+i})$ has an expansion
	\begin{equation}\label{EqOrderC}
		c=\vol\sum_ic_i\left[k^{n+1}+\left((n+1)i+\dfrac{\overline{S}}{2}\right)k^n\right]+O(k^{n-1})
	\end{equation}
	where $\overline{S}$ is the average of the scalar curvature of $\omega$.

	We conclude this section with a definition that will be useful for the estimates in the proof of Theorem~\ref{TheoMain2}.
	\begin{defin}[{\cite[p.~137]{rossthomas11b}}]\label{DefOrder}
		A sequence of real-valued functions $f_k$ on $X$ is of order $\Omega(k^p)$ if there exists a constant $\delta>0$ such that $f_k\geq\delta k^p$ uniformly on $X$ for all $p\gg 0$.
	\end{defin}

	\section{$C^\infty$-approximation of hermitian metrics}\label{sectionprooforbifold}
	As we anticipated in the Introduction, in this section we prove the approximation results in the \K setting needed in the proof of Theorem~\ref{TheoMain1}. 
	More precisely, given a cyclic \K orbifold $(X,\omega)$ polarized by $(L,h)$ with $2\pi\omega=i\de\deb\log h$, we get a $C^q$ approximation of $h$ (and consequently a $C^{q-2}$ approximation of $\omega$) by pulling back the Fubini-Study metric $h'_{FS}$ of Remark~\ref{RemDefRealFubini} via embeddings into a weighted projective space (cf. Proposition~\ref{PropConvRealFubini} below).
	Nevertheless, in order to highlight the connection to \cite{rossthomas11b}, we first prove the strengthening of their Theorem~4.6 to the $C^q$-norm.
	
	\begin{theor}\label{TheoMain2}
		Let $(X,L)$ be a polarized cyclic orbifold with $h$ a hermitian metric on $L$ whose curvature is $2\pi\omega$. 
		Denote with $Scal(\omega)$ the scalar curvature and with $\overline{S}$ its average. 
		Then for any given $q>1$ the pair $(h_{FS,k},\omega_{FS,k})$ $C^{q-2}$-converges to $(h,\omega)$ when $k$ tends to infinity. Namely,  we have
		\begin{equation}\label{EqApproxMetric}
			\dfrac{h_{FS,k}}{h}=1+\dfrac{\overline{S}-\Scal(\omega)}{2}k^{-2}+O(k^{-3})
		\end{equation}
		in the $C^q$-norm, and
		\begin{equation}\label{EqApproxForm}
			\omega=\omega_{FS,k}+O(k^{-2})
		\end{equation}
		in the $C^{q-2}$-norm.
	\end{theor}

	\begin{proof}[Proof of \eqref{EqApproxMetric}]
		Fix a $q\geq2$ and pick the coefficients $c_i$ in \eqref{EqCi} so that the expansion \eqref{EqExpansionCq} holds in $C^q$.

		We want to show that 
		\begin{equation}\label{EqExpansion}
			\alpha_k:=\dfrac{h_{FS,k}}{h}=1+\dfrac{\overline{S}-\Scal(\omega)}{2}k^{-2}+O\left(k^{-3}\right)\ \ \mbox{in}\ \ C^q.
		\end{equation}
		
		Notice that $\alpha_k$ satisfies the equation 
		\begin{equation}\label{EqImplicitAlpha}
			\sum_i(k+i)\alpha_k^{k+i}\BB=c
		\end{equation}
		where $c$ and $\BB$ are as above (cf. \cite[(4.10)]{rossthomas11b}).

		It was proven in \cite[page~138]{rossthomas11b} that $\alpha_k=1+O(k^{-2})$ and there exist positive constants $C_1$ and $C_2$ such that
		\begin{align}
			C_1\leq\alpha_k^j\ \  \mbox{ for all }&\frac{k}{2}\leq j\leq k;\\
			\alpha_k^j\leq C_2 \ \ \mbox{ for all }& 0\leq j\leq k.
		\end{align}
		Now define
		\begin{equation}\label{EqTauj}
			\beta_k:=1+\dfrac{\overline{S}-\Scal(\omega)}{2}k^{-2}+\sum_{j=1}^q \tau_jk^{-2-j}
		\end{equation}
		for smooth functions $\tau_j$ independent of $k$.
		Analogously to \cite[page~139]{rossthomas11b}, we get the equality
		\begin{align}\label{EqImplicitBeta}
			\nonumber    \sum_i(k+i)\beta_k^{k+i}\BB&=\sum_i(k+i)\left(1+\dfrac{\overline{S}-\Scal(\omega)}{2k}+O(k^{-2})\right)\BB\\
			\nonumber      &=vol\sum_i c_i\left[k^{n+1}+\left(\dfrac{\overline{S}-\Scal(\omega)}{2}+(n+1)i+\dfrac{\Scal(\omega)}{2}\right)k^n\right]+O(k^{n-1})\\
			\nonumber      &=vol\sum_i c_i\left[k^{n+1}+\left(\dfrac{\overline{S}}{2}+(n+1)i\right)k^n\right]+O(k^{n-1})\\
			&=c+O(k^{n-1})
		\end{align}
		in $C^0$ as $\BB=O(k^n)$ in $C^0$ by \eqref{EqDerivBergKer}.
		We now pick the functions $\tau_j$ in \eqref{EqTauj} so that the terms of order $n-j$ in \eqref{EqImplicitBeta} vanish for all $1\leq j\leq q$.
		In order to do so, notice that the coefficient of $k^{n-j}$ in \eqref{EqImplicitBeta} equals $b_0\tau_j+f_j$ where $f_j$ is independent of $k$ and of the functions $\tau_l$ for all $l\geq j$.
		Thus we can set $\tau_j=-\dfrac{f_j}{b_0}$ so that the terms of order $n-j$ in $k$ vanish for $1\leq j\leq q$. Hence we get
		\begin{equation}\label{EqImplicitBetaQ}
			\sum_i(k+i)\beta_k^{k+i}\BB=c+O(k^{n+p-q-1})\ \mbox{in} \ C^p\ \mbox{for}\ p=0,1,\ldots,q,
		\end{equation}
		where we have used \eqref{EqDerivBergKer}.
		Subtracting \eqref{EqImplicitBetaQ} from \eqref{EqImplicitAlpha} yields
		\begin{equation}\label{EqAlphaMinusBeta}
			(\alpha_k-\beta_k) \sum_i(k+i)\gamma_k\BB=O(k^{n+p-q-1})\ \mbox{in} \ C^p\ \mbox{for}\ p=0,1,\ldots,q,
		\end{equation}
		where
		\begin{equation}\label{EqGamma}
			\gamma_k=\sum_{j=1}^{k+i}\alpha_k^{j-1}\beta_k^{k+i-j}\, .
		\end{equation}
		Notice that $\gamma_k=\Omega(k)$ (cf. Definition~\ref{DefOrder}) since it is the sum of $k+i$ terms which are uniformly bounded from below. Thus the sum on the left hand side of \eqref{EqAlphaMinusBeta} is of order $\Omega(k^{n+2})$.
		In particular, for $p=0$ we get 
		
		\begin{equation}\label{EqC0Approx}
			\alpha_k-\beta_k=O(k^{-q-3})\ \mbox{in} \ C^0. 
		\end{equation}
		In fact, we claim that $D^p(\alpha_k-\beta_k)=O(k^{p-q-3})$ in $C^{p}$ for all $p=1,\ldots,q$, which yields \eqref{EqExpansion} for $p=q$.
		
		Our strategy is to prove the following facts by induction on $p=1,2,\ldots,q$, from which the claim follows 
		\begin{enumerate}[(i)]
			\item $D^p\alpha_k=O(k^{p-1})$;\label{Item1}
			\item $D^p\gamma_k=O(k^{p+1})$;\label{Item2}
			\item $D^{p}(\gamma_k\BB)=O(k^{n+p+1})$.\label{Item3}
			\item $D^p(\alpha_k-\beta_k)=O(k^{p-q-3})$ in $C^{p}$.\label{Item4}
		\end{enumerate}
		\textbf{Case $p=1$:} this is proved in \cite[page~140]{rossthomas11b}, but we recall the argument here for completeness.
		Differentiating \eqref{EqImplicitAlpha} yields
		$$D\alpha_k\sum_i(k+i)^2\alpha_k^{k+i-1}\BB=-\sum_i(k+i)\alpha_k^{k+i}D\BB\, .$$
		The sum on the left hand side is of order $\Omega(k^{n+2})$ while the sum on the right hand side is $O(k^{n+2})$ because $D\BB=O(k^{n+1})$ by \eqref{EqDerivBergKer}. Thus $D\alpha_k=O(1)$ which proves \eqref{Item1}. 
		
		Now for \eqref{Item2}  let us first notice that for $u,v\geq0$ with $u+v=k+i-1$ we have
		$$D(\alpha_k^u\beta_k^v)=u\alpha^{u-1}\beta_k^vD\alpha_k+v\alpha_k^u\beta_k^{v-1}D\beta_k.$$ 
		Therefore $D(\alpha_k^u\beta_k^v)=O(k)$ because $D(\alpha_k)=O(1)$, $D(\beta_k)=O(k^{-2})$ and both $\alpha_k$ and $\beta_k$ are uniformly bounded above.
		Now $D\gamma_k=O(k^2)$ follows from the fact that
		$$D\gamma_k=\sum_{j=1}^{k+i}D\left(\alpha_k^{j-1}\beta_k^{k+i-j}\right)  $$
		is a sum of $k+i$ terms of order $O(k)$, and this proves \eqref{Item2}.
		
		To show \eqref{Item3} we simply consider
		$$D(\gamma_k\BB)=D(\gamma_k)\BB+\gamma_kD(\BB)\, .$$
		Since $\gamma_k=O(k)$ and  $D(\BB)=O(k^{n+1})$, the claim follows from \eqref{Item2}.
		
		Now differentiating \eqref{EqAlphaMinusBeta} yields
		$$
		D(\alpha_k-\beta_k)\Omega(k^{n+2})=- (\alpha_k-\beta_k)\sum_i(k+i)D(\gamma_k\BB)+O(k^{n-q}) \ \mbox{in}\ C^1.$$
		Part \eqref{Item4} then follows from \eqref{EqC0Approx} and part \eqref{Item3}.

		\textbf{Inductive step:} Suppose that parts \eqref{Item1}-\eqref{Item4} hold for all natural numbers between $1$ and $p-1$. In particular, since part \eqref{Item4} holds, we get
		\begin{equation}\label{EqOrderDjAlpha}
			D^j(\alpha_k)=D^j(\beta_k)+O(k^{j-q-3})=O(k^{-2})\ \mbox{for all}\ 1\leq j\leq p-1.
		\end{equation}
		Let us now prove \eqref{Item1}. In order to do so, consider \eqref{EqImplicitAlpha} differentiated $p-1$ times, that is,
		\begin{equation}\label{EqDiffAlpha}
			D^{p-1}(\alpha_k)\sum_i(k+i)^2\alpha_k^{k+i-1}\BB=O(k^{n+p}).
		\end{equation}
		The right hand side is a sum of terms of order at most $O(k^{n+p})$. Each of these summands is a product involving some (or all) of the following factors
		\begin{enumerate}[(I)]
			\item $D^u\BB$ for $u\leq p-1$;\label{EqFactor1}
			\item powers of $D^v\alpha_k$ for $v\leq p-2$;\label{EqFactor2}
			\item powers of  $\alpha_k$;\label{EqFactor3}
			\item constants possibly involving powers of $k$.
		\end{enumerate}
		Now differentiating \eqref{EqDiffAlpha} we obtain the equation 
		\begin{align}\label{EqDiffAlpha2}
			\nonumber   D^{p}(\alpha_k)\Omega(k^{n+2})=&-D^{p-1}(\alpha_k)\sum_i(k+i)^2(k+i-1)\alpha_k^{k+i-2}D\alpha_k\BB\\
			-&D^{p-1}(\alpha_k)\sum_i(k+i)^2\alpha_k^{k+i-1}D\BB+D(RHS)
		\end{align}
		where $D(RHS)$ is obtained by differentiating the right hand side of \eqref{EqDiffAlpha}. Thus each term in $D(RHS)$ is obtained by differentiating one of the factors \eqref{EqFactor1}-\eqref{EqFactor3} in the summands of the right hand side of \eqref{EqDiffAlpha}.
		Notice that differentiating $\eqref{EqFactor1}$ raises the order of the summand by $1$ (again by \eqref{EqDerivBergKer}), differentiating $\eqref{EqFactor2}$ does not change the order of the summand and differentiating $\eqref{EqFactor3}$ lowers the order of the summand by $1$. We conclude that $D(RHS)$ has order at most $O(k^{n+p+1})$\footnote{In fact, the term of highest order in $k$ is $\sum_i(k+i)\alpha_k^{k+i}D^p\BB=O(k^{n+p+1})$.}.
		Since the order of the first two sums is $O(k^{n-1})$ and $O(k^{n+1})$ respectively, we have proven $\eqref{Item1}$.
		
		Now for \eqref{Item2} notice first that for $u,v\geq0$ with $u+v=k+i-1$ we have
		$$D^p(\alpha_k^u\beta_k^v)=u\alpha^{u-1}\beta_k^vD^p\alpha_k+O(k^{-2})=O(k^p)$$ 
		where the first equality follows from \eqref{EqOrderDjAlpha} and the second equality comes from part \eqref{Item1}.
		Thus \eqref{Item2} follows from the fact that
		$$D^p\gamma_k=\sum_{j=1}^{k+i}D^p\left(\alpha_k^{j-1}\beta_k^{k+i-j}\right)  $$
		is a sum of $k+i$ terms of order $O(k^p)$.
		
		Now also \eqref{Item3} follows from \eqref{Item2}, \eqref{EqDerivBergKer} and the inductive hypothesis  because 
		$$D^{p}(\gamma_k\BB)=\sum_{j=0}^{p}\binom{p}{j}D^{p-j}(\gamma_k)D^{j}(\BB)\, .$$
		In order to show \eqref{Item4} and conclude the proof, differentiate \eqref{EqAlphaMinusBeta} to get 
		\begin{equation*}
			\D^p(\alpha_k-\beta_k)\Omega(k^{n+2})= -\sum_{j=0}^{p-1}\binom{p}{j} \D^j(\alpha_k-\beta_k)\sum_i(k+i)\D^{p-j}(\gamma_k\BB)+O(k^{n+p-q-1}) \ \mbox{in}\ C^p.
		\end{equation*}
		Notice that $\D^{p-j}(\gamma_k\BB)=O(k^{n+p-j+1})$ and $\D^j(\alpha_k-\beta_k)=O(k^{j-q-3})$ so that we obtain \eqref{Item4} as wanted. 
	\end{proof}
	
	\begin{proof}[Proof of \eqref{EqApproxForm}]
		This proof is carried analogously to \cite{rossthomas11b} but we write it here for completeness.
		
		Applying $\de\deb\log$ to \eqref{EqApproxMetric} we get
		$$\omega_{h_{FS,k}}=\omega_h+O(k^{-2})=\omega+O(k^{-2})\ \mbox{in}\ C^{q-2}.$$
		Together with Lemma~\ref{LemFubiniStudy} this yields
		$$\omega_{FS,k}=\omega_{h_{FS,k}}+\dfrac{i}{2c}\de\deb f_k=\omega+\dfrac{i}{2c}\de\deb f_k+O(k^{-2})\ \mbox{in}\ C^{q-2}.$$
		Now, since $c$ is of order $\Omega(k^{n+1})$ by \eqref{EqOrderC}, it suffices to show that $f_k$ is constant up to terms of order $O(k^{n-1})$ in the $C^{q}$-norm.
		From the definition of the metric $h_{FS,k}$ it follows 
		$$f_k(x)=\vol\sum_i c_i\sum_\alpha\vert s^i_\alpha(x)\vert^2_{h_{FS,k}}=\vol\sum_i c_i \dfrac{h_{FS,k}^{k+i}}{h^{k+i}}\sum_\alpha\vert s^i_\alpha(x)\vert^2_{h}.$$
		Now by \eqref{EqApproxMetric} we get
		$$f_k(x)=\vol\sum_i c_i \left( 1+\dfrac{\Scal(\omega)-\overline{S}}{2k}+O(k^{-2})\right)\sum_\alpha\vert s^i_\alpha(x)\vert^2_{h}.$$
		Finally by \cite[Theorem~1.7]{rossthomas11a} we have the following $C^q$ expansion:
		$$\vol\sum_ic_i\sum_\alpha\vert s^i_\alpha(x)\vert^2_{h}=b_0k^n+b_1'k^{n-1}+\cdots$$
		where $b_0=\vol\sum_ic_i$ is a constant.
		Therefore we get
		$$f_k(x)=b_0k^n+O(k^{n-1}) \ \mbox{in}\ C^q$$ 
		which is what we were aiming for.
	\end{proof}
	Notice that Theorem~\ref{TheoMain2} is not the result we need for our purposes as $h_{FS}$ is not the hermitian metric whose curvature is proportional to the Fubini-Study form $\omega_{FS}$ on the weighted projective space.
	In order to apply this study to weighted Sasakian spheres we need to prove the convergence of the metrics $h'_{FS,k}$ obtained as pullback to $L$ of the metric $h'_{FS}$ defined in Remark~\ref{RemDefRealFubini}.
	
	\begin{prop}\label{PropConvRealFubini}
		Let $(X,L)$ be a polarized cyclic orbifold with $h$ a hermitian metric on $L$ whose curvature is $2\pi\omega$. 
		Then for a given $q>1$ the metrics $h'_{FS,k}$ converge to $h$ when $k$ tends to infinity. Namely, we have
		\begin{equation}\label{EqApproxRealFubini}
			\dfrac{h'_{FS,k}}{h}=1+O(k^{-1})
		\end{equation}
		in the $C^q$-norm.
	\end{prop}
	\begin{proof}
		Fix $q>1$ so that the coefficients $c_i$ in \eqref{EqCi} define a hermitian metric $h_{FS}$ as in \Cref{DefHermitianMetric}.
	By definition
		\begin{equation}
			h'_{FS}=e^{\frac{f}{2c}}h_{FS}
		\end{equation}
		
		where $f=\sum_i\sum_\alpha\vert t^i_\alpha\vert^2_{h_{FS}}$
		and  $\{t^i_\alpha\}$ is an orthonormal basis of $(V^*,\vert\cdot\vert_V)$.
		
		Therefore 
		\begin{equation}\label{EqDefRealFubini}
			h'_{FS,k}=e^{\frac{f_k}{2c}}h_{FS,k}
		\end{equation}
		where $f_k(x)=b_0k^n+O(k^{n-1})$ in $C^q$ as above.
		Moreover, since $c$ is of order $\Omega(k^{n+1})$ by \eqref{EqOrderC}, we have that $\dfrac{f_k}{2c}$ is a $O(k^{-1})$ in $C^q$.
		Therefore, $$ \dfrac{	h'_{FS,k}}{h_{FS,k}}=e^{\frac{f_k}{2c}}=1+O(k^{-1})\ \ \mbox{in}\ \ C^q.$$
		Finally, making use of \eqref{EqApproxMetric}, we have
		$$ \dfrac{	h'_{FS,k}}{h}=\dfrac{	h'_{FS,k}}{h_{FS,k}}\dfrac{	h_{FS,k}}{h}=\left(1+O(k^{-1})\right)\left(1+O(k^{-2})\right)=1+O(k^{-1})\ \ \mbox{in}\ \ C^q.$$
		The convergence of $\omega_{FS,k}$ to $\omega$ was already proved in Theorem~\ref{TheoMain2}.
	\end{proof}

	\section{Proof of Theorem~\ref{TheoMain1}}\label{sectionproof1}
	
	The line of argument to prove Theorem~\ref{TheoMain1} is the following.
	A Sasakian embedding of a Sasakian manifold $(M,\eta,R,g,\Phi)$ into a weighted Sasakian sphere determines an isometric \K embedding of its \K cone $(Y,J)$ into $\C^{n+1}\setminus\{0\}$. Our strategy is to define a sequence of holomorphic embeddings of the \K cone $(Y,J)$ into $\C^{n+1}\setminus\{0\}$ such that the induced \K cone metrics converge to the original one in the $C^q$-norm for a given $q$.
	This is done firstly for quasi-regular manifolds where we can exploit the fact that $Y$ is in fact a holomorphic line bundle over a \K orbifold and make use of Theorem~\ref{TheoMain2}. 
	
	Let $(M,\eta,R,g,\Phi)$ be a compact quasi-regular Sasakian manifold. By the discussion following Theorem~\ref{TheorStructure}, possibly after performing a transverse homothety, the \K cone of $M$ is the complement of the zero section of the dual $L^{*}$ of a orbiample line bundle $L$ over a compact \K orbifold $(X,\omega)$.
	In particular, the Sasakian structure on $M$ is determined by a hermitian metric $h$ on $L$ whose curvature is $2\pi\omega$ and $M$ identifies with the $U(1)$ orbibundle determined by $L^{*}$. As in Section~\ref{sectionbackground} denote with $\pi:M\lra X$ the Riemannian submersion given by the restriction to $M$ of the bundle projection $p:L^{*}\lra X$.

	By Theorem~\ref{TheoMain2} and Proposition~\ref{PropConvRealFubini} this, together with the choice of $q+2\in\N$, determines a sequence of embeddings $\phi_k$ into weighted projective spaces $(\PP(V),\omega_{FS})$ such that

	\begin{equation}\label{EqPullbackBundle}
		\phi_k^*\left(\OO_{\PP(V)}(1)\right)=L
	\end{equation}
	and   $\phi_k^*(h'_{FS},\omega_{FS})=(h'_{FS,k},\omega_{FS,k})$ converges to $(h,\omega)$ in $C^q$ as $k\lra\infty$.  
	Visually we have
	
	$$
	\begin{tikzcd}[column sep=large, row sep=large]
		(L, h'_{FS,k}) \arrow [r,"\widetilde{\phi}_k"] \arrow[d,"\pi_k"] &
		(\OO_{\PP(V)}(1),h'_{FS}) \arrow[d, "\pi_{FS}"] \\
		(X,\omega_{FS,k}) \arrow[r,"\phi_k"] & \left(\PP(V),\omega_{FS}\right)
	\end{tikzcd}
	$$
	where $\widetilde{\phi}_k$, i.e. the lift of the embedding $\phi_k$ to $L$, is an isometric \K embedding of \K cones.
	Moreover, the vertical maps restricted to the $U(1)$-orbibundles are Boothby-Wang fibrations.
	Now the hermitian metric $h'_{FS,k}$ defines a Sasakian structure $(\eta_k,R_k,g_k,\Phi_k)$ on $M$ as in \eqref{EqFromHermitianToSasakian}. 
	Furthermore, these structures converge to $(\eta,g)$ in the $C^q$-norm because the metrics $h'_{FS,k}$ converge to $h$ in $C^{q+2}$ and so do the associated coordinates $t_k$ on the \K cone $Y=L^*\setminus\{0\}$.
	Notice that we could pull back the \K cone structure on $\OO_{\PP(V)}(-1)=V\setminus\{0\}$ after performing a $\D$-homothetic transformation so that the induced Sasakian structure do converge to the original metric on $M$.
	This concludes the proof in the case where the Sasakian structure $(M,\eta,g)$ is quasi-regular.
	
	Suppose now that the structure $(M,\eta,R,g,\Phi)$ is irregular. By Proposition~\ref{PropRukimbira}, any irregular structure can be $C^\infty$-approximated by a sequence of quasi-regular structures $(M,\eta_j,R_j,g_j,\Phi_j)$ obtained by type-I deformation. 
	Moreover, by the previous part each of these admits a $C^q$-approximation by structures $(M,\eta_{j,k},R_{j,k},g_{j,k},\Phi_{j,k})$ induced by embeddings $\varphi_{j,k}$ into a weighted Sasakian sphere.
	Therefore, the Sasakian structures induced on $M$ by the embeddings $\varphi_k:=\varphi_{k,k}$ into weighted spheres converge $C^q$ to the structure $(\eta,g)$. This concludes the proof of Theorem~\ref{TheoMain1}.


\begin{thebibliography}{99}
		
	
		
		
		
		\bibitem{bandecappellettimontanoloi20}
		G.~Bande, B.~Cappelletti--Montano and A.~Loi,
		\emph{{$\eta$}-{E}instein {S}asakian immersions in non-compact {S}asakian space forms},
		Ann. Mat. Pura Appl. (4) {\bf 199} (2020), no. 6, 2117--2124. 
		
		
		
		
		
		
		
		\bibitem{boyergalicki08}
		C.~P.~Boyer and K.~Galicki,
		{\sl Sasakian geometry},
		Oxford Mathematical Monographs. Oxford University Press, Oxford, 2008.
		
		
		
		
		
		
		\bibitem{cappellettimontanoloi19}
		B.~Cappelletti-Montano and A.~Loi,
		\emph{Einstein and {$\eta $}-{E}instein {S}asakian submanifolds in spheres},
		Ann. Mat. Pura Appl. (4) {\bf 198} (2019), no. 6, 2195--2205.
		
		
		
		\bibitem{collins18}
		T.~C.~Collins and G.~Sz\'ekelyhidi,
		\emph{{K}-semistability for irregular Sasakian manifolds},
		J. Differential Geom. (1) {\bf 109} (2018) 81--109.
		
		
		
		\bibitem{DaiLiuMa06}
		X.~Dai, K.~Liu and X.~Ma,
		\emph{On the asymptotic expansion of {B}ergman kernel},
		J. Differential Geom. {\bf 72} (2006), no. 1, 1--41.
		
		
		
		
		\bibitem{futakionowang09}
		A.~Futaki, H.~Ono and G.~Wang,
		\emph{Transverse {K}\"{a}hler geometry of {S}asaki manifolds and toric	{S}asaki-{E}instein manifolds},
		J. Differential Geom. {\bf 83} (2009), no. 3, 585--635. 
		
		
		
	
		
		
		\bibitem{hesun16}
		W.~He and  R.~Sun,
		\emph{Frankel conjecture and Sasaki geometry},
		Adv. Math. {\bf 291} (2016), 912--960. 
		
		
		\bibitem{HerrmannHsiaoLi17}
		H.~Herrmann, C.-Y.~Hsiao and  X.~Li,
		\emph{Szeg\H{o} kernel expansion and equivariant embedding of {CR}
			manifolds with circle action},
		Ann. Global Anal. Geom. {\bf 52} (2017), no. 3, 313--340.
		
		
		
		\bibitem{HerrmannHsiaoLi22}
		H.~Herrmann, C.-Y.~Hsiao and  X.~Li,
		\emph{Szeg\H{o} kernels and equivariant embedding theorems for {CR} manifolds},
	Math. Res. Lett. {\bf 29} (2022), no. 1, 193--246.
	
		
		
		\bibitem{HerrmannHsiaoMarinescuShen22}
		H.~Herrmann, C.-Y.~Hsiao, G.~Marinescu and W.-C.~Shen,
		\emph{Semi-classical spectral asymptotics of Toeplitz operators on CR manifolds},
		Preprint arXiv:2303.17319v2.
		
		
		
		\bibitem{kotschickplacini22}
		B.~Kotschick and G.~Placini,
		\emph{Sasaki structures distinguished by their basic Hodge numbers},  
		Bull. London Math. Soc. https://doi.org/10.1112/blms.12667
		
		
		
		
		
		
		
		\bibitem{marinescuma}
		G.~Marinescu and X.~Ma,
		{\sl Holomorphic Morse inequalities and Bergman kernels},
		Progress in Mathematics, Series Volume 254, 2007, Birkh\"auser Basel.
		
		
		
		
		\bibitem{martellisparksyau08}
		D.~Martelli, J.~Sparks and  ST.~Yau,
		\emph{Sasaki–Einstein Manifolds and Volume Minimisation},
		Commun. Math. Phys. {\bf 280} (2008), 611--673.
		
		
		
		\bibitem{orneaverbitsky05}
		L.~Ornea and  M.~Verbitsky,
		\emph{An immersion theorem for Vaisman manifolds},
		Math. Ann. {\bf 332}, (2005) 121--143.
		
		\bibitem{orneaverbitsky07a}
		L.~Ornea and  M.~Verbitsky,
		\emph{Embeddings of compact Sasakian manifolds},
		Math. Res. Lett. {\bf 14}  no. 4, (2007), 703--710. 
		
		\bibitem{orneaverbitsky07b}
		L.~Ornea and  M.~Verbitsky,
		\emph{Sasakian structures on CR-manifolds},
		Geom. Dedicata {\bf 125} (2007), 159--173. 
		
		\bibitem{orneaverbitsky22}
		L.~Ornea and  M.~Verbitsky,
		\emph{Principles of Locally Conformally K\"ahler Geometry},
		Preprint arXiv:2208.07188
		
		\bibitem{placini21}
		G.~Placini,
		\emph{Sasakian immersions of Sasaki-Ricci solitons into Sasakian space forms},
		Journal of Geometry and Physics \textbf{166}, 104265 (2021).
		
		
		
		\bibitem{placini22}
		G.~Placini,
		\emph{Approximation of regular Sasakian manifolds},
		Pacific J. Math. {\bf 327} (2023), no.1, 167--181.
		
		
		
	
		
		\bibitem{rossthomas11a}
		J.~Ross and R.~Thomas
		\emph{Weighted {B}ergman kernels on orbifolds},
		J. Differential Geom. {\bf 88} (2011), no. 1, 87--107. 
		
		
		\bibitem{rossthomas11b}
		J.~Ross and R.~Thomas
		\emph{Weighted projective embeddings, stability of orbifolds, and constant scalar curvature {K}\"{a}hler metrics},
		J. Differential Geom. {\bf 88} (2011), no. 1, 109--159. 
		
		
		\bibitem{ruan98}
		W.~D.~Ruan
		\emph{Canonical coordinates and Bergmann metrics},
		Comm. Anal. Geom. {\bf 6}  no. 3, (1998), 589--631. 
		
		
		\bibitem{rukimbira95}
		P.~Rukimbira
		\emph{Chern–Hamilton’s conjecture and K-contactness},
		Houston J. Math. {\bf 21} (1995), no. 4, 709--718.
		
		
		\bibitem{sparks11}
		J.~Sparks
		{ \sl Sasaki-Einstein manifolds},
		Surveys in differential geometry. Volume XVI. Geometry of special holonomy and related topics, 265--324, Int. Press, Somerville, MA, 2011. 
		
	
		
		\bibitem{tian90}
		G.~Tian,
		\emph{On a set of polarized K\"ahler metrics on algebraic manifolds},
		J. Differential Geom. {\bf 32} (1990), 99--130.
		
		
		
		\bibitem{zelditch98}
		S.~Zelditch,
		\emph{Szeg\"o kernels and a theorem of Tian},
		International Mathematics Research Notices, Volume 1998, {\bf 6}, (1998), 317--331.
		
		
		
		
		
		
	\end{thebibliography}
\end{document}